
\input phyzzx
\hfuzz 50pt

\font\mybb=msbm10 at 12pt
\def\bb#1{\hbox{\mybb#1}}

\def\bC {\bb{C}}
\def\bR {\bb{R}}

\def\bQ {\bb{Q}}

\def\bI {\bb{I}}
\def\bJ {\bb{J}}

\def\bfomeg{\omega\kern-7.0pt\omega}
\def\bfOmeg{\Omega\kern-8.0pt\Omega}


\REF\polb{J. Polchinski, {\it String Theory}, Vol I and II,
 Cambridge University Press (1998).}
\REF\towna{P.K. Townsend, {\it The Eleven-Dimensional 
 Supermembrane Revisited }, Phys. Lett. 
 {\bf B350} (1995) 184; hep-th/9501068.}
 \REF\wittena{E. Witten, {\it  String 
 Theory Dynamics in Various Dimensions},
 Nucl. Phys. {\bf B443} (1995) 85; hep-th/9503124.}
 \REF\banks{
 T. Banks, W. Fischler, S.H. Shenker and L. Susskind, 
 {\it M Theory as a Matrix Model: A Conjecture},
Phys.Rev. {\bf D55} (1997) 511; hep-th/9610043.}   
 \REF\hulltown{C.M. Hull and P.K. Townsend, {\it  Unity 
 of Superstring Dualities},
  Nucl. Phys. {\bf B438} (1995) 109; hep-th/9410167.} 
  \REF\gibbonsc{ G.W. Gibbons, G. T.
   Horowitz and P.K. Townsend,
  {\it Higher Dimensional Resolution of 
  Dilatonic Black Hole Singularities},
Class. Quant. Grav. (1995) 297; hep-th/9410073.} 
\REF\gppkt{G. Papadopoulos and P.K. Townsend,  
{\it Intersecting M-branes}, Phys. Lett.
 {\bf B380} (1996) 273; hep-th/9603087.}
\REF\duff{
M.J. Duff and K.S. Stelle, {\it Multi-membrane 
solutions of d = 11 supergravity},
Phys.Lett. {\bf B253} (1991) 113.}
\REF\guven{R. G\"uven, {\it Black p-brane solutions 
of D=11 supergravity theory},
Phys. Lett. {\bf 276B} (1992) 49.}
\REF\cjs{E. Cremmer, B. Julia
 and J. Scherk, {\it Supergravity
 Theory in 11 
Dimensions}, Phys. Lett. {\bf 76B} (1978) 409.}
\REF\chs{
C. G. Callan, Jr., J. A. Harvey and A. Strominger, 
{\it  Worldbrane actions 
for string solitons},
 Nucl.Phys. {\bf B367} (1991)60.}
 \REF\calic{J. Darok, R. Harvey and 
 F. Morgan, {\it Calibrations on $\bR^8$},
Tranc. Am. Math. Soc., Vol. 307 {\bf 1} (1988), 1.}
\REF\ptb{G. Papadopoulos and A. Teschendorff, 
{\it Multi-angle five-brane intersections},
Phys. Lett. {\bf B443} (1998) 159; hep-th/9806191.}
\REF\ptc{G. Papadopoulos and A. Teschendorff, 
{\sl Grassmannians, Calibrations and 
Five-Brane Intersections}; hep-th/9811034.} 
\REF\paptown{G. Papadopoulos and P.K. Townsend,
{\it Compactification of D = 11 
Supergravity on Spaces of Exceptional
Holonomy},  Phys.Lett. {\bf B357} (1995) 300; 
 hep-th/9506150.}
 \REF\papad{G. Papadopoulos, {\it  Brane 
 Surgery with Thom Classes},
 JHEP (1999) 9905:020; hep-th/9905073.}  
\REF\tseytlinb{A.A. Tseytlin, {\it Harmonic 
Superposition of M-Branes}, Nucl. Phys.{\bf B475}
(1996) 149; hep-th/9604035.}
\REF\twist{P.S. Howe and 
G. Papadopoulos, {\it Twistor Spaces
for HKT Manifolds}, Phys. Lett. 
{\bf B379} (1996) 80; hep-th/9602108.}
\REF\rocek{S.J. Gates, Jr., C.M. Hull and M. Ro\v cek,
{\it Twisted Multiplets and New
 Supersymmetric Nonlinear Sigma Models},
 Nucl. Phys. {\bf B248} (1984) 157.} 
 \REF\howea {P.S. Howe and 
 G. Papadopoulos, {\it Ultraviolet
 Behaviour of Two-Dimensional
Nonlinear Sigma Models},
 Nucl. Phys. {\bf B289} (1987) 264 \&
{\it Further Remarks on the
Geometry of Two-Dimensional
 Nonlinear Sigma Models}, Class. 
Quantum Grav. {\bf 5} (1988)
1647.}
\REF\coles{ R. A. Coles and
 G. Papadopoulos, {\it The Geometry 
of the One-dimensional
Supersymmetric Non-linear Sigma Models}, 
Class. Quantum Grav. {\bf 7} (1990)
427-438.}
\REF\gibbonsp{G.W. Gibbons, 
G. Papadopoulos  and K.S. Stelle, 
{\it HKT and OKT Geometries on 
Soliton Black Hole Moduli Spaces},
 Nucl.Phys. {\bf B508} (1997)623; hep-th/9706207.}
\REF\hitchin{ N.J. Hitchin, 
A. Karlhede, U. Lindstr\"om and M. Ro\v cek,
{\sl  Hyperkahler Metrics and Supersymmetry}, 
 Commun.Math.Phys. {\bf 108} (1987) 535.}
 \REF\opferhowe{P.S. Howe, 
 A. Opfermann and G. Papadopoulos, 
 {\it Twistor Spaces for QKT Manifolds},
 Commun.Math.Phys. {\bf 197}(1998)713; hep-th/9710072.}
 \REF\ivanovd{S. Ivanov, {\it Geometry 
 of Quaternionic K\"ahler 
 connections with torsion},  ICTP Preprint (1999) No. IC99195.}
 \REF\salamon{S.M. Salamon, {\it Quaternionic 
 K\"ahler Manifolds} 
 Invent. Math. {\bf 67} (1982) 143 \&
 {\it Differential Geometry of 
 Quaternionic Manifolds}, Ann. Sci.
 Ecole Norm. Sup (4) {\bf 19} (1986) 31 \& 
  {\it Riemannian Geometry 
  and Holonomy Groups },
 Pitman Res. Notes in Math. 
 201 Longman, Harlow 1989.}
 \REF\ivanov{B. Alexandrov and S. Ivanov, {\it Vanishing 
 Theorems on Hermitian
Manifolds}, math/9901090.}
\REF\ivanova{B. Alexandrov,
 G. Grantcharov  and S. Ivanov, {\it 
The Dolbeault operator on 
hermitian surfaces}; math/9902005.}
\REF\proeyen{P. Spindel, A. Sevrin, 
W. Troost and A. Van Proeyen, 
{\it Extended Supersymmetric Sigma 
Models on Group Manifolds. 1. The
Complex Structures}
Nucl. Phys. {\bf B308} (1988) 662. }
\REF\opferman{
A. Opfermann and G. Papadopoulos, {\it 
 Homogeneous HKT and QKT Manifolds}; math-ph/9807026.}
\REF\joyce{D. Joyce, {\it Compact 
Hypercomplex and Quaternionic
Manifolds}, J. Diff. Geom. {\bf 35} (1992) 743.}
\REF\poon{G. Grantcharov 
and Yat-Sun Poon, {\it Geometry of Hyper-K\"ahler 
Connections with Torsion}; math.dg/9908015. } 
\REF\micha{J. Michelson and 
A. Strominger, {\it Superconformal 
Multi-Black Hole
Quantum Mechanics}, JHEP (1999) 
9909:005; hep-th/9908044.}
\REF\gutpap{J. Gutowski and 
G. Papadopoulos, {\it The dynamics 
of very special black holes}, Phys. Lett. {\bf B} in press; 
hep-th/9910022 \&
{\it Moduli spaces for four- and five-
dimensional black holes}; hep-th/0002242.}
\REF\cali {R. Harvey and H. Blaine Lawson, Jr,
{\it Calibrated Geometries}, Acta Math. {\bf
148} (1982) 47.}
\REF\ptg{J. P. Gauntlett, G.W. Gibbons, 
G. Papadopoulos and P.K. Townsend, {\sl
Hyper-K\"ahler Manifolds and Multiply 
Intersecting Branes}, Nucl. Phys. {\bf B500} (1997)
133; hep-th/9702202.}
\REF\jerome{J.P. Gauntlett, D.A. Kastor and J. 
Traschen, {\it Overlapping Branes in M-theory}
 Nucl. Phys.  {\bf B478} (1996) 544; hep-th/9604179.}
 \REF\yau{A. Strominger,
  S-T. Yau and E. Zaslow, {\it Mirror 
Symmetry is T-Duality},
 Nucl. Phys.{\bf B479} (1996) 243;
 hep-th/9606040.}
\REF\papas{G.W. Gibbons 
and G. Papadopoulos, {\it Calibrations
 and Intersecting Branes}, Commun.
Math. Phys. {\bf 202} (1999) 593; hep-th/9803163.} 
\REF\west{J.P. 
Gauntlett, N.D. Lambert and
 P.C. West, {\it Branes and Calibrated 
Geometries}, Commun.Math.Phys. {\bf 202} (1999) 571; 
hep-th/9803216.} 
\REF\spence{A. Acharya, 
J.M. Figueroa-O' Farrill and B. Spence, {\sl 
Branes at Angles and 
Calibrated Geometry}, JHEP (1998) 9804:012;  
hep-th/9803260 .}


\Pubnum{ \vbox{ \hbox{} } }
\date{ March 2000}
\pubtype{}

\titlepage

\title {\bf Brane Solitons and Hypercomplex Structures}

\author{G. PAPADOPOULOS}
\address{Department of Mathematics,
\break King's College London,
\break
Strand,
\break
London WC2R 2LS, U.K.}

\abstract{The investigation of strings and M-theory 
involves the understanding of various
 BPS solitons  which in a certain
approximation can be thought of as
 solutions of ten- and eleven-dimensional
supergravity theories.  These solitons
 have a brane or a intersecting brane
interpretation, saturate a  bound and are associated with
parallel spinors with respect to a connection
 of the spin bundle  
of spacetime.  A class of intersecting 
brane configurations is examined and it is
shown that the geometry of 
spacetime is hyper-K\"ahler with torsion.
 A relation between
these hyper-K\"ahler geometries with 
torsion and  quaternionic calibrations is 
also demonstrated. }

\endpage
\pagenumber=2

\chapter{Introduction}

The main achievement in theoretical 
physics the past few years is the realization
that all five string theories [\polb] 
are related amongst themselves and that are limits
of a another theory which has been 
called M-theory [\towna, \wittena].
The precise nature of M-theory  
remains a mystery but by now an impressive amount of
evidence has been gather which point to its existence.
Most of these arise from
 investigating the low energy approximation
of strings and M-theory which are 
described by ten- and eleven-dimensional
supergravities, respectively; see however [\banks].
The use of supergravity theories 
in this context is two-fold. First,
the conjectured  duality symmetries of string 
theories are discrete subgroups 
of the continuous duality symmetries 
of the field equations of supergravity 
theories [\hulltown]. Second, the supergravity
theories admit solutions which have the 
interpretation of extended objects
embedded in the spacetime, called branes, which are  the
BPS solitons of strings and M-theory.
 A consequence  of their BPS property
is that branes are  stable under deformations of the various
parameters of the theories, 
like for example coupling constants. Because
of this, they can be used to 
compare the various limits of M-theory and thus
establish the relations amongst 
the various string theories.
Extrapolating from the properties 
of eleven-dimensional and ten-dimensional
 supergravities, M-theory
is thought to have the following essential ingredients:

 \item{\bullet} A low energy 
 description in terms of the eleven-dimensional
supergravity,

\item{\bullet} M-2- and M-5-branes, and

\item{\bullet} limits that describe 
all five string theories in ten dimensions.

The  solutions of  supergravity theories with a brane 
interpretation  have some common properties. 
The associate spacetime of a p-brane
has an asymptotic region isometric
 to $\bR^{(1,p)}\times \bR^n$,  where,
 viewing the p-brane as a (p+1)-dimensional submanifold of
spacetime, $\bR^{(1,p)}$ and 
$\bR^n$ are identified with the worldvolume
and transverse directions 
of the p-brane, respectively;
 $(p+n)$ is equal either to nine (string theory) 
or to ten (M-theory). In the above  
asymptotic region a mass $m$  and a charge $q$ 
 per unit $\bR^p\subset \bR^{(1,p)}$ volume is
defined, i.e $m$ and $q$ are energy 
and charge densities, respectively. 
Then using the properties of 
supergravity theory, a bound  
can be established [\gibbonsc] as
$$
m\geq \alpha\, |q|\ ,
$$
where for string theory branes 
$\alpha$ depends on the string coupling constant 
$\lambda$. The precise dependence 
of $\alpha$ distinguishes between the various
types of branes
as follows: $\alpha\sim \lambda^0$ for 
fundamental strings, $\alpha\sim \lambda^{-1}$
for Dirichlet branes or D-branes for 
short and  $\alpha\sim \lambda^{-2}$ for
Neveu-Schwarz 5-branes or  NS-5-branes for short.
The solutions that are of interest are 
those that saturate the above bound
leading to BPS type  configurations. BPS configurations
 are associated with parallel spinors 
 with respect to a connection
which occurs naturally in 
supergravity theories. The BPS branes of strings
and M-theory admit sixteen parallel spinors.
Apart from the stability of
these BPS solutions that has 
already been mentioned above, superposition rules
have been found that allow to 
combine two or more such solutions
and construct new ones [\gppkt].
 The solutions that arise from superpositions
of BPS branes also admit parallel 
spinors which typically are less than
those of the branes involved in the superposition.

In this paper,  the geometry of a class
of BPS brane solutions of 
supergravity theories and that of
their superpositions will be described. 
I shall begin with a 
description of BPS M-2-brane 
[\duff] and M-5-brane [\guven] solutions
 of eleven-dimensional supergravity [\cjs].
Then I shall explain the connection 
between BPS solutions and parallel spinors.
I shall continue with the NS-5-brane 
solution of type II ten-dimensional
supergravity theories [\chs] and show 
that the geometry of this solution
is  hyper-K\"ahler with torsion (HKT).
Then I shall explore the various superpositions of
 NS-5-branes and I shall demonstrate 
 that these superpositions
are  related to the quaternionic 
calibrations in $\bR^8$ [\calic].
 I shall interpret these superpositions 
 as intersecting NS-5-branes
 and I shall show that the geometry 
 of these solutions is again HKT. 
 Most of these results have 
 appeared in [\ptb , \ptc ]. Finally,
 I shall state my conclusions.

\chapter{Eleven-Dimensional Supergravity}

I shall not attempt to give a full 
description of eleven-dimensional
supergravity. This can be found in 
the original paper of Cremmer, Julia and
Scherk who constructed the theory [\cjs].
 Here I shall only emphasize some
aspects of the geometric structure of the theory.
In field theoretic terms, the theory describes
the dynamics of the graviton $g$, 
a three-form gauge potential $A$ and
a gravitino $\psi$. The latter is a
 spinor-valued one-form which does
not enter in the analysis below and 
so it will be neglected in what follows.
Geometrically, let $(N; g, F, \nabla)$
 be an eleven-dimensional
 spin manifold $N$ of  signature 
 $(-,+,\dots,+)$ equipped 
with a metric $g$, a closed 
four-form $F$, locally $F=dA$, 
and a connection $\nabla$.
In the physics literature $\nabla$ 
is called superconnection and
$$
\nabla\, : \qquad C^\infty(S)
\rightarrow \Omega^1(N)\otimes C^\infty(S)\ ,
$$
where $S$ is the spin bundle over $N$ 
and ${\rm rank}\,(S)=32$.
This connection can be written as
$$
\nabla=D+T(F)\ ,
$$
$D$ is the connection of $S$ induced from 
the Levi-Civita connection of
the metric $g$ and
$$
T_M(F)dx^M =-{1\over 144}F_{NPQR} 
\big(\Gamma_M{}^{NPQR}-8\delta^N{}_M \Gamma^{PQR}\big)
dx^M\ ,
$$
where $\{\Gamma^M; M=0,\dots, 10\}$ is a 
basis in the Clifford algebra
${\rm Clif}(1,10)$ and 
$\Gamma^{M_1M_2\cdots M_n}=
\Gamma^{[M_1} \Gamma^{M_2}\cdots \Gamma^{M_n]}$.
The dynamics of the theory is 
described\foot{The conventions for forms are as
follows:
$\omega={1\over p!} \omega_{a_1, \dots, a_p}
 dx^{a_1}\wedge\dots\wedge dx^{a_p}$,
$|\omega|^2={1\over p!} 
\omega_{a_1, \dots, a_p} \omega^{a_1, \dots, a_p}$
and ${}^*{}^*\omega=-(-1)^{p(n-p)}\omega$, 
where ${}^*$ is the
Hodge star operation and
 $n=11$ in the present case.}
 by the action 
$$
S=\int\, d^{11}x\,\sqrt{|{\rm det}\,g|}\big(R(g)-
2 |F|^2\big)-{4\over 3}A\wedge F\wedge F\ ,
$$
where $R(g)$ is the Ricci scalar 
of the metric $g$ and the norm of $F$
is taken with respect to $g$.
The above action consists from the 
 Einstein-Hilbert term, the standard
kinetic term for $F$ and a Chern-Simons term, respectively.
The equations of the fields $g$ 
and $A$ can be derived by varying the above
action.

There are two classes of 
solutions to the field equations
depending on whether or not 
$F=0$. If $F=0$, then the field 
equations imply that the Ricci 
tensor of $g$ vanishes. Therefore a 
large class of solutions is $\bR^{(1,10-n)}\times M^n$, where 
$M^n$ is a manifold of
 appropriate special holonomy, i.e $SU(k)$, 
$k=2,3,4$ $(n=2k)$; $Sp(2)$ $(n=8)$; $G_2$ $(n=7)$; Spin(7) 
$(n=8)$.  Such solutions 
admit parallel spinors and
 have found application in the various 
compactifications of M-theory [\paptown]. 
The other class of solutions is that for which $F\not=0$.
For such solutions to
 have a brane interpretation, 
it is required that they have 
an asymptotic region which is isometric to either 
$\bR^{(1,2)}\times \bR^8$ 
or $\bR^{(1,5)}\times \bR^5$. The former 
asymptotic behaviour is that of 
 M-2-brane while the latter 
 is that of M-5-brane. Then after 
imposing appropriate decaying
 conditions on the fields as they 
approach these asymptotic
 regions,  the charges per unit volume
of the M-2- and M-5-branes 
can be defined   as follows: 
$$
q_2=\int_{S^7}\,({}^*F-A\wedge F)\ ,
$$
where $S^7\subset \bR^8$, and
$$
q_5=\int_{S^4}\, F\ , 
$$
where $S^4\subset \bR^5$, respectively. 
Then adapting the positive mass theorem
of general relativity to this case,
 the bounds can be established, 
$$
m_2\geq\alpha_2 |q_2|
$$
or
$$
m_5\geq\alpha_5 |q_5|\ ,
$$ 
where $m_2$ and $m_5$ are the 
M-2-brane and M-5-brane masses per unit 
volume, respectively.
The manifolds that saturate those 
bounds admit sixteen parallel spinors
with respect to the connection $\nabla$.

To be specific, the BPS 
M-2-brane solution [\duff] is
$$
\eqalign{
ds^2&= h^{{1\over3}} \big(h^{-1} 
ds^2(\bR^{(1,2)})+ds^2(\bR^8)\big)
\cr
{}^*F&=\mp{1\over2}\star dh\ ,}
$$
where 
$$
h=1+{q_2\over |y|^6}
$$
is a harmonic function on $\bR^8$, 
$y\in \bR^8$, the Hodge star operation on $F$
is with respect to the metric $g$ 
on $N$ and the Hodge star operation on $dh$ is 
with respect to the flat 
metric on $\bR^8$. The M-2-brane is located
at $y=0$.
There are two asymptotic regions. 
One is as $|y|\rightarrow \infty$,
where the spacetime $N$ becomes 
isometric to $\bR^{(1,2)}\times \bR^8$ as expected.
The other is as $|y|\rightarrow 0$
 in which case $N$ becomes isometric to
$AdS_4\times S^7$; $AdS_4$ is a 
Minkowski signature analogue of the
standard hyperbolic four-manifold. 
It turns out though that
the BPS membrane solution
 develops a singularity behind the
$AdS_4\times S^7$ region.

The BPS M-5-brane solution [\guven] is
$$
\eqalign{
ds^2&= h^{{2\over3}}
 \big(h^{-1} ds^2(\bR^{(1,2)})+ds^2(\bR^5)\big)
\cr
F&=\mp{1\over2}\star dh\ ,}
$$
where 
$$
h=1+{q_5\over |y|^3}
$$
is a harmonic function on 
$\bR^5$, $y\in \bR^5$ and the Hodge
star operation on $dh$ is 
with respect to the flat metric on $\bR^5$. 
The M-5-brane is located
at $y=0$.
There are also two asymptotic regions in this case. 
One is as $|y|\rightarrow \infty$
where the spacetime $N$ becomes isometric to 
$\bR^{(1,5)}\times \bR^5$ as expected.
The other is as $|y|\rightarrow 0$ 
in which case $N$ becomes isometric to
$AdS_7\times S^4$; $AdS_7$ is a 
Minkowski signature analogue of the
standard hyperbolic seven-manifold. 
The BPS M-5-brane solution is not
singular.

Despite much progress in constructing 
solutions of eleven-dimensional
supergravity which admit parallel
 spinors, there is no systematic way
to tackle the problem of constructing solutions
of the theory with any number of 
parallel spinors. Of course 
this is related to understanding the
properties of the connection $\nabla$. 
To give an example of the subtleties
involved, let us consider the case of
 solutions with thirty two parallel spinors
which is the maximal number possible.
 A straightforward example of
a such spacetime is the Minkowski 
$N=\bR^{(1,10)}$ space with vanishing $F$.
However, this is not all. There are
 two more cases, $N=AdS_4\times S^7$
with $F$ the volume form of $AdS_4$
 and $N=AdS_7\times S^4$ with $F$ the
volume form of $S^4$.

A related and still unresolved problem
 is to construct {\sl localized} solutions
of eleven-dimensional supergravity 
theory with the interpretation of intersecting
M-2- and M-5-branes. In particular 
consistency of the M-theory picture
suggests that there should be a solution 
that has the interpretation of
a M-2-brane ending orthogonally on a 
M-5-brane associated with
eight parallel spinors. No such solution 
has been found so far; For a recent
discussion of this see [\papad]. 
Other BPS solutions involving M-5-branes 
and M-2-branes have been found though,
like a solution which has the 
interpretation of a M-2-brane \lq passing through'
a M-5-brane [\tseytlinb].

\chapter{Type IIA Strings}

A  geometric insight into the 
properties of the connection
$\nabla$ of eleven-dimensional 
supergravity can be given after reducing
  M-theory to type IIA string theory.
It turns out that in a sector of 
the supergravity theory
associated with type IIA strings,
 IIA supergravity, the connection $\nabla$
of the spin bundle is   induced 
from certain connections
with torsion
of the tangent bundle of spacetime.

For this consider the 
eleven-dimensional spacetime
$N=S^1\times M$ with
$$
\eqalign{
ds^2(N)&= e^{{4\over3}\phi} d\theta^2+
 e^{-{2\over3}\phi} ds^2(M)
\cr
F&=d\theta\wedge H\ ,}
\eqn\red
$$
where $\theta$ is the angle which 
parameterizes the circle of radius 
$r$. It is assumed that the vector
 field $X={\partial\over\partial \theta}$
is an isometry of $N$ which 
leaves in addition $F$ invariant. 
In field theoretic terms, the
 metric $\gamma$ in $ds^2(M)$ describes
the graviton in ten dimensions,
 the closed three-form $H$ is the
$NS\otimes NS$ form field 
strength and $\phi$ is the dilaton;
$\gamma, H$ and $\phi$ are the so called
common sector fields of string theory.
  The IIA supergravity has 
  additional fields
which can also be derived from 
eleven dimensions but the above
sector will suffice for our 
purpose. The type IIA string coupling constant
 is related 
to the radius of the circle $S^1$ as
$\lambda=r^{3\over2}$ [\wittena]. 
So for small radius, the  
string coupling constant is small
and M-theory reduces to IIA strings.

The dynamics of the common sector 
fields in ten dimensions
 can be described by the action
$$
S=\int\, d^{10}x\, e^{-2\phi}
\,{\sqrt{|{\rm det}\,\gamma}|}\, \big( R(\gamma)-2
 |H|^2+4 |d\phi|^2\big)\ ,
$$
where the norms are taken with respect to the
 metric $\gamma$ on $M$.

As it has been mentioned, a 
simplification occurs in the
description of the connection 
$\nabla$ using $N=S^1\times M$ and \red\ as above.
For this, first observe that 
the spin bundle $S$  decomposes
as $S=S^+\otimes S^-$, where 
$S^+$ and $S^-$ are spin bundles
over $M$ with 
${\rm rank }\,S^+={\rm rank }\,S^-=16$. This
is due to the decomposition of 
the  spinor representation of ${\rm spin}(1,10)$
into the sum of the two 
irreducible spinor representations
of ${\rm spin}(1,9)$. Next it 
turns out that the connection $\nabla$ decomposes
into two connections one on 
$S^+$ and one on $S^-$ which are induced by
the connections
$$
\nabla^\pm=D\pm H
$$
of the tangent bundle, respectively,
where $D$ is the Levi-Civita 
connection of the metric $\gamma$.
 The connections $\nabla^\pm$ are metric
connections with torsion the
 closed three-forms $\pm H$, respectively.
There are two more conditions 
that arise from reducing the connection
$\nabla$ of eleven-dimensional 
supergravity to IIA supergravity which
involve the dilaton $\phi$. 
However, these two conditions do not
give rise to additional 
restrictions on the parallel spinors of $\nabla^\pm$
connections in the examples
 that we shall investigate below. So we shall
not consider them further.

The above simplification in the
 structure of the connection 
$\nabla$ has some profound
consequences. One of them  is 
that the existence of parallel 
spinors with respect to the
 connection $\nabla$ of $S$
  depends
on the holonomy of the 
connections $\nabla^\pm$ of the tangent bundle of $M$.

The NS-5-brane solution of IIA supergravity [\chs] is
$$
\eqalign{
ds^2(M)&=ds^2\big(\bR^{(1,5)}\big)+ h ds^2({\bQ})
\cr
H&=\mp{1\over2}\star dh
\cr
e^{2\phi}&=h\ ,}
$$
where the Hodge star 
operation on $dh$ is with respect
to the flat metric on $\bR^4$ 
and $h$ is a harmonic function 
on $\bR^4=\bQ$, $\bQ$ is the 
quaternionic line, 
$$
h=1+{1\over |q|^2}\ ,
$$
$q\in \bQ$. The NS-5-brane is located at $q=0$.
The spacetime $M$ is diffeomorphic 
to  $\bR^{(1,5)}\times \big(\bQ-\{0\}\big)$
and it has two asymptotic regions, 
$\bR^{(1,5)}\times \bR^4$ as 
$|q|\rightarrow \infty$ and
 $\bR^{(1,5)}\times \bR\times S^3$ as
$|q|\rightarrow 0$. In what follows we 
shall choose $H=-{1\over2}\star dh$.

The non-trivial part of the metric 
of $M$ is that on $\bQ-\{0\}$.
To investigate the geometry on 
 $\big(\bQ-\{0\}\big)$, we introduce two
hypercomplex structures 
${\bf I}=\{I_1, I_2, I_3\}$ 
and ${\bf J}=\{J_1, J_2, J_3\}$ as follows:
$$
\eqalign{
I_1(dq)&=i\, dq\ , 
\qquad I_2(dq)=j\, dq\ ,\qquad I_3(dq)=k\, dq\ ,
\cr
J_1(dq)&=- dq\, i\ , 
\qquad J_2(dq)=- dq\, j\ ,\qquad J_3(dq)=- dq\, k\ ,}
$$
where $i,j,k$ are the imaginary
 unit quaternions. Observe
that the two hypercomplex
 structures commute, $[{\bf I}, {\bf J}]=0$.
The metric on $\bQ-\{0\}$ is
 hermitian with respect to both
 hypercomplex structures. In
  addition, the hypercomplex structure
 ${\bf I}$ is compatible with the $\nabla^-$ connection 
 ($\nabla^-{\bf I}=0$) and the 
 hypercomplex structure ${\bf J}$ is 
 compatible with the $\nabla^+$ connection 
 ($\nabla^+{\bf J}=0$); Note that 
 the torsion $H$ has support on $\bQ-\{0\}$.
 Therefore the holonomy of $\nabla^\pm$
 is in $SU(2)$. In fact the holonomy of  $\nabla^\pm$ 
 is $SU(2)$ and so the NS-5-brane 
 admits sixteen parallel spinors.
 This fact follows from
  representation theory. As it will be demonstrated
 shortly, 
the geometry of the NS-5-brane
 can be summarized
by saying that it admits two 
commuting strong HKT structures.

\chapter{Hyper-K\"ahler Manifolds with Torsion}

Let $(M, g, {\bf J})$ be a 
Riemannian hyper-complex manifold with metric
$g$ and hypercomplex structure 
${\bf J}$; ${\rm dim} M=4k$.
The manifold $(M, g, {\bf J})$ 
admits a HKT structure [\twist] if

\item{\bullet} The metric $g$ is
 hermitian with respect to all three
complex structures.

\item{\bullet} There is a compatible 
connection $\nabla$ with both the metric $g$
and the hypercomplex structure ${\bf J}$
 which has torsion a {\sl three form} $H$.

There are two types of HKT structures 
on manifolds, the {\sl strong} and the 
{\sl weak}, depending
on whether or not the torsion 
three-form is a closed, respectively.

Torsion has appeared in the
 physics literature since the
early attempts to 
incorporate it in a relativistic
theory of gravity. In supersymmetry, 
metric connections with torsion 
a closed three-form have appeared  
in the context  of IIA and
 IIB supergravities but the
relation to HKT geometry was not established.
Connections with torsion a {\sl closed} 
three-form  appeared next  in the
investigation of 
two-dimensional supersymmetric 
sigma models [\rocek, \howea]. For a 
 class of models, the sigma model 
 manifold satisfies  conditions
which can be organized in 
one or two copies of what it is 
now called  {\sl strong} HKT structure.
The general case of connections 
with torsion {\sl any} three-form 
were found in the investigation
of  one-dimensional supersymmetric sigma models
 [\coles, \gibbonsp]. For a 
class of models, the sigma model manifold
satisfies conditions which can be organized 
 as one or two
copies of what it is now called
{\sl weak} HKT structure. 
The definition of the HKT
 structure as a new structure
on manifolds was given in
[\twist]. In the same paper, 
 the  strong and weak 
HKT structures were introduced, a formulation of a
 HKT structure in terms of 
conditions on K\"ahler forms was given,
and a twistor construction for 
the HKT manifolds was proposed.
The latter two properties of  
HKT manifolds are similar to those
 of   hyper-K\"ahler  
  manifolds [\hitchin]. There is also
 a generalization of the Quanternionic 
 K\"ahler structure on manifolds to include
 torsion. The Quaternionic K\"ahler
  manifolds with torsion (QKT) have
 been introduced in  [\opferhowe] 
 and further investigated
 in [\ivanovd]. QKT manifolds 
 admit a twistor construction similar
 to that of Quaternionic K\"ahler ones [\salamon].

A straightforward consequence  of the  
definition of HKT manifolds is that
the holonomy of the connection $\nabla$ 
is in $Sp(k)$. Some other developments
related to these connections with 
three-form torsion are the vanishing 
theorems of [\ivanov, \ivanova] for
 certain cohomology groups. 
Many examples of HKT manifolds
 have been constructed.
These include a class of group 
manifolds with  strong HKT structures in [\proeyen].
Specifically, the Hopf surface 
$S^1\times S^3$ admits two strong HKT structures.
Homogeneous weak HKT manifolds 
have been constructed in [\opferman] using
the hypercomplex structures of 
[\joyce]. Inhomogeneous weak HKT structures have
been given
on $S^1\times S^{4k-1}$ in [\poon].

In physics, the NS-5-brane solution constructed in 
the previous section clearly admits two strong
HKT structures each associated with the  
hypercomplex structures 
on $\bQ-\{0\}$ defined by left and
 right quaternionic multiplication, respectively.
Other examples are certain 
(strong and weak) HKT structures that appear on the
moduli spaces of five-dimensional 
black holes [\gibbonsp, \micha, \gutpap].

The HKT structure has many properties 
some of them found in [\twist]
and more have been derived in [\poon]. 
One of them is the following:
Let $M$ be  hypercomplex manifold 
with respect to ${\bf J}$ and
equipped with a three-form $H$. 
$M$ admits a HKT structure if
$$
d\omega_{{\bf J}}-2i_{{\bf J}} H=0\ ,
\eqn\lem
$$
where $\omega_{{\bf J}}$ are the 
three K\"ahler forms associated with
the hyper-complex structure and 
$i_{{\bf J}}$ are the inner derivations
with respect to the three complex 
structures. This equation will
be used later to construct
 new HKT manifolds. In fact if two of the
above conditions are satisfied,
 they imply the third.
Observe also that the torsion 
of an HKT manifold can be specified
from the metric and the complex
 structures [\twist]. So in what follows,
we shall not give the expression for the torsion.  

\chapter {Quaternionic Calibrations}

Calibrations have been 
introduced by Harvey and Lawson [\cali] 
to construct a large class
of minimal submanifolds. Here, I shall use  
calibrations to find a new class
of solutions of IIA 
supergravity that 
has the interpretation 
of intersecting branes.
This new class of solutions 
admits a strong HKT structure.

A calibration of degree k is a k-form $\omega$ such 
that for every k-plane $\eta$
in $\bR^n$
$$
\omega({\buildrel\rightarrow \over {\eta}})\leq 1\ ,
$$ 
where ${\buildrel\rightarrow \over {\eta}}$ 
is the co-volume  form
of $\eta$.

The contact set $G_\omega$ of a 
calibration is the subset of ${\rm Gr}(k, \bR^n)$
of k-planes that saturate the 
above bound. Usually $G_\omega$ 
is a homogeneous space.
There are many examples of calibrations, like K\"ahler and
Special Lagrangian,  which have been 
extensively investigated both in 
 mathematics and physics.
In the present case, the relevant
 class of calibrations are
the quaternionic calibrations
that have been described by Daroc, 
Harvey and Morgan in [\calic].
For this, we identify $\bR^8=\bQ^2$ 
and introduce the 
hypercomplex structures $\bI=\{I_1, I_2,
I_3\}$ and
$\bJ=\{J_1, J_2, J_3\}$ on $\bQ^2$
induced by left and right quaternionic
 multiplication, respectively, i.e.  
$$
\eqalign{
I_1(du)&=i\,du\ ,
\qquad I_2(du)=j\,du\ ,\qquad I_3(du)=k\, du\ ,
\cr
J_1(du)&=-du\, i\ ,
\qquad J_2(du)=-du\, j\ ,\qquad J_3(du)=-du\, k\ ,}
$$
where $i,j,k$ are the 
imaginary unit quaternions.
In addition, we define
$$
\phi_{\bJ}={1\over3} 
\sum^3_{r=1} \omega_{{}_{J_r}}\wedge \omega_{{}_{J_r}}\ ,
$$
where $\omega_{{}_{J_r}}$ is 
the K\"ahler form of the $J_r$ complex structure
with respect to the standard 
metric on $\bQ^2$, and similarly $\phi_{\bI}$ for
$\bI$. 
There are several calibration 
forms that can be constructed 
using the K\"ahler
forms above. The calibration 
forms of quaternionic
calibrations are found by 
averaging the calibration form
$\phi_{\bJ}$ with various 
calibration forms constructed using the
$\bI$ hypercomplex structure.
 These calibration forms and 
 their contact sets are 
 summarized in the following table:

\leftline{~}
\centerline{{\underbar {Quaternionic 
Calibrations in $\bR^8$}}}
$$
\vbox{\settabs 2\columns 
\+{\rm Calibration Form} $\omega$& 
{\rm Contact Set} $G_\omega$  \cr
\+${1\over2} \big(\phi_{\bI}+\phi_{\bJ}\big)$
&\qquad$S^1$\cr
\+${1\over5} \omega_{I_1}\wedge\omega_{I_1}
 +{3\over5}\phi_{\bJ}$&\qquad$S^2$\cr
\+${1\over4}\Omega+{3\over4}\phi_{\bJ}$&\qquad $S^3$\cr
\+$\phi_{\bJ}$&\qquad $S^4$\cr
}
$$
\leftline{~}

In the second case above, 
instead of $I_1$ any other of the 
complex structures of $\bI$
can be used.
In the third case, the form 
$\Omega$ is the $Spin(7)$ invariant form
associated with $\bI$. In the 
last case, the contact set is
the grassmannian of quaternionic 
lines in $\bQ^2$, $Gr(1; \bQ^2)=S^4$.
 The contact sets of the 
 quaternionic calibrations
are computed by observing
 that the groups that leave invariant
the above forms act transitively 
on the calibrated planes. For more details
see [\calic]. 

\chapter{HKT Geometries and Branes}

The strategy of constructing HKT geometries
 in eight dimensions is to superpose
a  HKT geometry on $\bQ-\{0\}$ along  
four-planes
in $\bR^8$. The  four-planes which will 
be used are in the contact sets of 
quaternionic calibrations that has been
 described in the previous section.
The construction [\ptb, \ptc]  involves 
the following steps:

\item{(i)} Consider the HKT metric
$$
d\sigma^2={1\over |q|^2} dq d\bar q
\eqn\thkt
$$
 on $\bQ-\{0\}$, where $q\in \bQ$. The 
 torsion of this HKT geometry
is that of the NS-5-brane that we have described.

\item{(ii)} Introduce the maps
$$
\eqalign{
\tau:\qquad &\bQ^2\rightarrow \bQ
\cr
&u\rightarrow \tau(u)=p_1u^1+p_2u^2-a\ ,}
$$
where $p_1, p_2, a$ are quaternions.

\item{(iii)} Define the metric 
$$
ds^2=\sum_{\tau} r^2(\tau) \tau^*d\sigma^2+ds^2(\bQ^2)\ ,
$$
where the sum is over a finite number of maps $\tau$, 
$r(\tau)\in \bR$ and $ds^2(\bQ^2)$
is the standard flat metric on $\bQ^2$.

The manifold 
$K=\bQ^2-\cup_{\tau}\{\tau^{-1}(0)\}$ 
admits a HKT structure
associated with the hypercomplex
 structure ${\bf J}$ induced
 by the hypercomplex structure 
 $\bJ$ on $\bQ^2$. To show this, the torsion of 
the HKT structure on $K$ is given
by pulling back the torsion of the
 NS-5-brane with respect to
the maps $\tau$ and then summing 
up over the various
maps $\tau$.  The key observation 
is that the differential
$d\tau$ commutes with the
 hypercomplex structures on $\bQ-\{0\}$
and $K$ defined by {\sl right}
 quaternionic multiplication.  Using
this, the condition \lem\ for
 $K$ can be written as 
$$
\sum_\tau\tau^*\big(d\omega_{{\bf J}}-2i_{{\bf J}} H\big)
=0\ ,
$$
where the expression inside the 
brackets is that  of the
condition \lem\ for the HKT
 structure (i) on $\bQ-\{0\}$ and so vanishes
identically. Thus $K$ admits a
 HKT structure with respect to the
$\nabla^+$ connection and 
${\bf J}$ hypercomplex structure.

Observe that $d\tau$ for generic parameters
$\{p_1, p_2\}$ does not commute with the hypercomplex
structure induced by left quaternionic 
multiplication on $\bQ-\{0\}$ and $K$.
So $K$ does not admit a HKT structure
 with respect to this hypercomplex
structure.

To make connection with the calibrations 
of the previous sections as promised,
we consider a HKT geometry constructed 
using several maps $\tau$ with generic
parameters $(p_1, p_2;a)$.
Such a HKT geometry is independent
from the parameterization of the maps $\tau$ 
and depends only on the arrangements
of quaternionic lines $\tau^{-1}(0)$ in $\bQ^2$. 
To see this, observe that two  maps $\tau$ and $\tau'$ 
give the same HKT structure if their
parameters are related as
$$
(p_1', p_2'; a')= (up_1, up_2; ua)
$$
for some $u\in \bQ$, $u\not=0$. So the 
inequivalent HKT structures 
associated with each map $\tau$ are 
parameterized by the bundle space
of the canonical quaternionic line bundle over the 
Grassmannian ${\rm Gr}(1; \bQ^2)$. 
In turn the quaternionic lines 
defined by ${\rm Ker}\,d\tau$ are
in ${\rm Gr}(1; \bQ^2)$ which
 is precisely the contact set of
last calibration in the table
 given in the previous section. Observe
that the calibration form and 
the HKT connection $\nabla^+$ are compatible
with the same hypercomplex structure ${\bf J}$. 

The HKT geometries that we are 
considering are complete provided that
the subspaces $\tau^{-1}(0)$ are
 in general position. Near the
intersection of two such subspaces, 
the HKT metric is isometric to that
of $(\bQ-\{0\})\times (\bQ-\{0\})$,
 where the metric on $\bQ-\{0\}$
is given as in \thkt.

Finally, the above HKT geometries can be used to
construct new  solutions of IIA supergravity  as
$$
\eqalign{
ds^2(M)&=ds^2(\bR^{(1,1)})+ds^2(K)
\cr
e^{2\phi}&= ({\rm det}\gamma)^{{1\over4}}\ .}
$$ 
The IIA supergravity three-form field 
strength $H$ is given is terms of
the torsion of the HKT manifold $K$.
The brane interpretation of such 
solution is that of NS-5-branes intersecting
on a string. The NS-5-brane 
associated with the map $\tau$ is
 located at $\tau^{-1}(0)$.

\chapter{Special Cases}

As we have seen for a generic choice of 
maps $\tau$ the  HKT geometries in eight 
dimensions found in the
previous section were associated
with quaternionic lines in $\bQ^2$ 
given by ${\rm Ker}\, d\tau$. This 
establishes a correspondence
between HKT geometries and quaternionic 
calibrations with calibration form
$\phi_{\bJ}$. This correspondence
 can be extended to the rest
of the quaternionic calibrations.
 For this, instead
of considering generic maps $\tau$ with parameters
$(p_1, p_2;a)$ to construct the HKT geometries,
we restrict them in an appropriate 
way. There are four cases to consider,
including the HKT geometry on $K$ that 
has been mentioned above,
as illustrated in the following table:

\leftline{~}
\centerline{{\underbar {HKT Geometries 
in Eight Dimensions}}}
$$
\vbox{\settabs 2\columns 
\+{ $\bar p_1 p_2$}& {\rm Ker}~ $d\tau$   \cr
\+$\bR$&\quad$S^1$\cr
\+$\bC$&\quad$S^2$\cr
\+{\rm Im}$\bQ$&\quad $S^3$\cr
\+$\bQ$&\quad $S^4$\cr
}
$$
In the first column we denote the
 restriction on the parameters of the
map and in the second column the 
set that ${\rm ker}~ d\tau$ lies as we
vary the map $\tau$ in the same 
class. Comparing
the above table with that which 
contains the contact sets
of  quaternionic calibrations in 
section five, we observe that  ${\rm ker}\,d\tau$
lies in a contact set in all four cases. 

The holonomy groups of the connections
$\nabla^\pm$ in each of the above cases
 are given in the following table:
$$
\vbox{\settabs 5\columns 
\+{ $\bar p_1 p_2$:}& $\bR$&$\bC$&${\rm Im}\bQ$&$\bQ$   \cr
\+$\nabla^-$:&$Sp(2)$&$SU(4)$&$Spin(7)$&$SO(8)$\cr
\+$\nabla^+$:&$Sp(2)$&$Sp(2)$&$Sp(2)$&$Sp(2)$\cr
}
$$

The holonomy of $\nabla^+$ is in $Sp(2)$ 
in all cases because it is compatible
with the ${\bf J}$ hypercomplex structure. 
Now if $\bar p_1 p_2$ is real for
all $\tau$ involved in the construction of
 HKT geometry, then $d\tau$
commutes with the hypercomplex structures 
of $K$ and $\bQ-\{0\}$ induced
by {\sl left} quaternionic multiplication. 
This leads to another HKT structure
on $M$ compatible with the $\nabla^-$ 
connection. So the holonomy
of $\nabla^-$ is in $Sp(2)$ as well. 
This HKT geometry was 
found in [\ptg] and a special case 
in [\jerome]. If $\bar p_1 p_2$ is 
complex, say, with respect
to the $I_1$ complex structure, then $d\tau$
commutes with the complex 
 structures of $K$ and $\bQ-\{0\}$ induced
by left quaternionic multiplication 
with the quaternionic unit $i$.
This makes the $\nabla^-$ 
connection compatible with the $I_1$
complex structure which implies
 that the holonomy of $\nabla^-$
is in $U(4)$. In fact it turns out 
that the holonomy
of $\nabla^-$ is in $SU(4)$. 
Observe that the K\"ahler form of $I_1$ appears in
the construction of the calibration form in this case.
A similar analysis can be 
done for the remaining case.

Some of the above HKT  geometries can
 be related to toric hyper-K\"ahler
geometries [\hitchin, \ptg]. In particular, the HKT 
geometries associated with maps $\tau$ 
such that $\bar p_1 p_2\in \bR$ are  
 T-dual (mirror symmetry)
to toric hyper-K\"ahler manifolds 
[\gibbonsp].  In this case mirror symmetry
transforms manifolds of one class, 
HKT manifolds, to manifolds of
another class, hyper-K\"ahler 
manifolds. This is because the T-duality
above is performed as many times 
as the number of tri-holomorphic vector
fields of the toric hyper-K\"ahler
 manifold which is {\sl less} than the middle
dimension of the manifold. This is 
unlike the case of mirror symmetry
for Calabi-Yau spaces where T-duality
 is performed in as many directions as
the middle dimension of the manifold [\yau].
The HKT geometries with 
$\bar p_1 p_2\in \bR$ also appear in the context 
moduli spaces of a class 
of black holes in five dimensions [\gibbonsp].

It is of interest to ask the question
 whether it is possible to construct
supergravity solutions which have
 the interpretation
of intersecting branes using other 
calibrations from those employed
above. To be more specific
instead of the quaternionic calibrations,
 one may also use K\"ahler
or special Lagrangian calibrations to do the 
superposition.
Unfortunately, in both these  cases, a superposition
similar to that employed for 
 quaternionic calibrations does not lead 
to solutions of supergravity
 field equations. This may due
to the fact that the resulting geometries
 depend on the particular
parameterization of the maps $\tau$.
On the other hand string perturbation theory 
considerations seem
to suggest that superpositions of
 the kind employed above lead to BPS
brane configurations [\papas, \west, \spence].
However, it is not known how to
 construct in a systematic way the
corresponding  supergravity
solution from a BPS brane 
configuration of string theory.

\chapter{Concluding Remarks}

The understanding of the 
non-perturbative properties
of string theory requires the 
investigation of various solitons.
In the low energy approximation 
these solitons have an interpretation
as branes or as intersecting
 branes and are solutions of
various supergravity theories. 
A class of such solutions was 
presented and their
construction was related to 
quaternionic calibrations.

The problem of finding the 
intersecting brane solutions
of supergravity theories 
 has not been tackled in complete generality.
Although many examples of 
such solutions are known, there
does not seem to be a
 systematic way to find a solution
for each BPS brane configuration of
string theory. The 
resolution of this will require a
better understanding of the 
supercovariant derivative 
of supergravity theories.
The method of  calibrations 
that I have presented 
led  to the construction of a large class
of these solutions but it 
has  limitations some of which  has already
mentioned. However, the 
solutions that we are seeking for which
the form field strengths 
do not vanish ($F\not=0$)
are in the same universality 
class as hyper-K\"ahler, Calabi-Yau and
other special holonomy 
manifolds as far as the holonomy
of the supercovariant 
connections of the supergravity theories is
concerned. So it may be 
that powerful methods of 
algebraic and differential
geometry that have been 
developed to construct examples of the
latter may be extended 
to find examples of the former.

After the end of the conference,  
the  moduli space of a  class 
of five-dimensional
black holes was determined
 and it was found to be a
weak HKT manifold [\micha]. This result 
was further generalized 
in [\gutpap] to a larger class 
of four- and five-dimensional
black holes.  The moduli spaces 
of all five-dimensional
black holes that admit at least four 
parallel spinors are 
 HKT manifolds.

\leftline{~}
\leftline{~}
\leftline{~}
\noindent{\bf Acknowledgments:}   I would 
like to thank the organizers
of the Second Meeting  on \lq\lq Quaternionic 
Structures in Mathematics
and Physics" for their kind invitation 
and their warm welcome at the conference. 
This work is partly
 supported by the PPARC grant
   PPA/G/S/1998/00613. G.P. is 
supported by a University 
Research Fellowship
from the Royal Society. 

\refout
\bye